\documentclass[a4paper, 12pt, 
twoside, 
onecolumn, 
nochapterprefix, %
normalheadings, 
halfparskip, 
headsepline, 
1.0headlines, 
headinclude, %
footinclude, %
DIV15, 
liststotoc, 
bibtotoc,
cleardoubleempty] 
{scrartcl} 

\usepackage[a4paper,left=25mm,right=25mm,top=20mm,bottom=15mm]{geometry}
\usepackage[manualmark,headsepline,plainheadsepline]{scrpage2}
\usepackage{floatflt,multicol}
\usepackage{amsmath, amssymb, amsthm, latexsym, eurosym, ngerman}
\usepackage{graphicx, tikz, array}
\usepackage{makeidx, eurosym}
\usepackage[latin1]{inputenc}
\usepackage{tikz,eurosym}

\begin{document}

\clearscrheadfoot 
\ihead{\headmark}
\ohead[\pagemark]{\pagemark}
\pagestyle{empty}

\title{Zum 200. Geburtstag \\ von Evariste Galois}
\author{Markus Wessler}
\maketitle

\begin{multicols}{2}
Es gibt Menschen, die gleichzeitig innerhalb wie auch außerhalb ihrer Zeit leben, und nicht selten ist gerade das Werk eines solchen Lebens, einer solchen Lebensleistung, besonders nachhaltig. Evariste Galois, der vor 200 Jahren, am 25.10.1811, geboren wurde, hat in seinem sehr kurzen Leben nachhaltig gewirkt, und er war ein Mensch seiner Zeit: Geboren zwischen zwei Revolutionen -- von denen die Ideen der einen, der ersten, großen, stark auf ihn wirkten und in deren Geist er erzogen wurde, und von denen die zweite seine beiden letzten Lebensjahre prägte -- war er ein sehr politischer Mensch. Und doch stand er auch stets ein wenig wie ein Zuschauer am Rande seines eigenen Lebens.  

Nun teilt Galois definitiv das Schicksal vieler großer Mathematikerinnen und Mathematiker, die der Breite der Öffentlichkeit eher weniger bekannt sind. Und 2011 wird sicher eher als Kleist-Jahr oder (schon wieder) Mahler-Jahr, weniger aber als Galois-Jahr wahrgenommen. Man könnte sagen, Galois habe die Algebra revolutioniert: ein großes Wort. Kratzt man jedoch all das Romantische und Romantisierende, das Mystische und Mystifizierende ab, das sich im Lauf der Jahre um den Kern der mathematischen Leistung Galois' gelegt hat, so bleibt tatsächlich erstaunlich viel Substanz übrig. Die Galoistheorie ist heute ein Pfeiler der modernen Mathematik. Studierende der Mathematik begegnen ihr in der Regel im dritten oder vierten Semester. Die Mathematik wird oft als \glqq eindeutig\grqq\ oder \glqq unfehlbar\grqq\ bezeichnet. Vielen mag nicht bewusst sein, dass sie aber auch eine lange Geschichte mit vielen Wandlungen, auch Irrungen und Wirrungen, hinter sich hat. Dieser Wandlungsprozess wird vor allem offenkundig, wenn wir die Mathematik als Sprache betrachten. Archimedes, Pythagoras, Thales haben anders über Mathematik gesprochen, ihre Erkenntnisse auf andere Weise festgehalten, als Kopernikus und Kepler oder Leibniz und Newton oder Gauß und Euler dies taten -- und erst recht anders, als Galois dies dann getan hat. 

Das Jahrhundert vor Galois, das 18. Jahrhundert: eine Phase, übervoll von mathematischer Aktivität. Es wurden ungeahnte Strukturen und Zusammenhänge entdeckt. Galois' Hauptverdienst ist es, zu dieser strukturierenden Strömung einen erheblichen Teil beigetragen zu haben. Er brachte das uralte Feld der algebraischen Gleichungen in Zusammenhang mit der Theorie der Gruppen: ein Begriff, der in seiner heutigen Bedeutung überhaupt erst von Galois geprägt wurde. Galois schuf aber nicht nur Neues oder Neusprachliches, sondern durch seine Theorie erschienen klassische Probleme der Geometrie in einem völlig neuen Licht und konnten elegant gelöst werden: so etwa die Würfelverdoppelung oder die Winkeldreiteilung. Die strukturschaffenden, strukturerhellenden Verbindungen zwischen Algebra und Geometrie: ein enormer Schritt auf dem Weg der modernen Mathematik.

Evariste Galois wird am 25.10.1811 in Bourg-La-Reine, wenige Kilometer südlich von Paris, geboren und wächst im behüteten Umfeld eines philosophisch, literarisch, religiös gebildeten Bürgertums auf. In seinem Elternhaus ist eine Schule untergebracht, die sein Großvater vor der Revolution gegründet hat und die nun sein Vater, Nicolas-Gabriel Galois leitet. Dieser genießt unter seinen Mitbürgern ein enormes Ansehen und wird im Jahr 1815, während der 100 Tage von Napoleon, zum Bürgermeister gewählt, ein Amt, das er überraschenderweise auch unter König Ludwig behalten darf. 

Der kleine Evariste vergnügt sich häufig mit dem Verfassen von Reimen und Dialogen, die er bei Familienfesten, zusammen mit seiner älteren Schwester Nathalie und seinem jüngeren Bruder Alfred, darbietet. Die außerordentliche mathematische Begabung lässt sich aus den wenigen verlässlichen Quellen dieser Zeit allerdings noch nicht herauslesen. Ab seinem 12. Lebensjahr besucht er das Lyceum \glqq Louis-le-Grand\grqq\ in Paris, so wie viele große Mathematiker Frankreichs vor ihm: Borel, Hadamard, Hermite, Lebesgue. Evariste tritt nur zögerlich in deren Fußstapfen; am Ende seines zweiten Jahres aber erhält er vier Auszeichnungen. Bald schon begehrt er auf, wird selbstbewusster, zeigt erste Anzeichen eines unbequemen, liberalen Geistes: rebelliert gegen den Druck von oben. Kann er sich (heftig unterstützt durch einen Brief seines Vaters an den Direktor Laborie) noch im Sommer 1826 gegen eine Rückstufung wehren, so muss er sich im Januar 1827 Labories Druck beugen und ab sofort die \glqq Seconde\grqq\ besuchen. Doch ist genau dies der schicksalhafte Moment in seinem jungen Leben: Er hat nun die Möglichkeit, erstmals an den damals noch nicht regulär angebotenen Mathematikkursen teilzunehmen, und er erlebt gewissermaßen eine Initialzündung, als er dort das Buch \glqq Elements de Géometrie\grqq\ von Legendre in die Hände bekommt und verschlingt. 

Der Einzug, den die Mathematik hier in Galois' Welt hält, stellt einen spürbaren Wendepunkt in seinem jungen Leben dar. Schon zuvor wenig kommunikativ, zieht er sich nun beinahe völlig in selbstauferlegte Isolation zurück, um sich ganz den neuen Ideen zu widmen, die rasch von seiner ganzen Vorstellung Besitz ergreifen. Er lernt die neu entdeckte Sprache schnell, formalisiert sie gar weiter und verfasst erste eigene Arbeiten in einem so zuvor nicht gekannten, knapp gehaltenen Stil. Es ist auch und gerade dieser Stil, der ihm später Probleme bereiten soll. Mit kaum 16 Jahren widmet er sich fortan in der Hauptsache einem der klassischen algebraischen Probleme: der Auflösbarkeit von Gleichungen fünften (und höheren) Grades durch elementare Formeln. Je mehr er sich in die Welt der Algebra hineinversenkt, um so mehr kommt er seinem Umfeld abhanden. 

Galois entwickelt nun im Prinzip schon die Grundlagen der nach ihm benannten Theorie. Bereits bekannt ist zu diesem Zeitpunkt, dass es für eine allgemeine Gleichung fünften Grades keine elementare Lösungsformel gibt; dies hat der Norweger Niels-Henrik Abel, den Galois sehr schätzt, einige Jahre zuvor gezeigt. Galois untersucht nun, wie sich die Nullstellen einer allgemeinen Gleichung zueinander verhalten, wenn sie permutiert werden. Später entdeckt und benennt er erstaunliche Gesetzmäßigkeiten und entwickelt schließlich den Begriff der \glqq auflösbaren Gruppe\grqq\, den er der Auflösbarkeit von Gleichungen gegenüberstellt.

Die Ecole Polytechnique ist auch damals schon eine der angesehensten elitären Hochschulen Frankreichs; vor allem das Ingenieurstudium gilt als ein Eckpfeiler. Galois nimmt hier im Juni 1828 an der Aufnahmeprüfung teil -- und fällt durch. Dies wird oft als Ursprung alles späteren Unheils angenommen. Man kann dieser Meinung sein, aber dieser Schluss ist nicht zwingend; hier ist viel Spekulation.  Bei vielen Menschen mag eine schwere Enttäuschung in frühen Jahren, in ihrer noch nicht so großen Erfahrungswelt, einen großen Raum einnehmen -- und ihr Leben endet doch nicht zwangsläufig in einer Katastrophe. Galois jedenfalls geht sein Leben nun erst recht zielstrebig an und plant eine Wiederholung der Prüfung im folgenden Jahr. Dazu besucht er Kurse bei dem begnadeten Mathematiklehrer Louis Richard, der von Galois' Fortschritten in der Mathematik schnell begeistert ist und ihn bei seiner ersten Veröffentlichung unterstützt: einer Arbeit über stetige Funktionen in den \emph{Annales de Gergonne}. Zwei weitere Arbeiten von Galois überreicht Richard direkt dem Leiter der Akademie der Wissenschaften, Louis Augustin Cauchy. 

Richard ist davon überzeugt, dass die Ecole Polytechnique schließlich gar nicht anders könne, als Galois aufzunehmen -- worin  er nicht Recht behält. Kurz nachdem Galois mit dem Selbstmord des Vaters einen schweren Schlag hinzunehmen hat, fällt er im Sommer 1829 erneut durch die Aufnahmeprüfung.  Über den Verlauf dieser Prüfung gibt es unterschiedliche Darstellungen, mehr oder weniger dramatische. Allen Beschreibungen gemein ist aber die Grundhaltung Galois': Von sich selber in höchstem Maße überzeugt, kommt er mit den Fragen der Prüfer nicht zurecht; sie sind ihm zu banal, er verweigert die Antwort entweder komplett oder lässt sich zu einigen lediglich verbalen Erklärungen -- den Gebrauch von Schwamm und Kreide lehnt er ab -- in seinem knappen Stil herab. 

Das Jahr 1830 wird für Galois nicht besser als das Jahr 1829: Im Februar wird er zwar an der Ecole Préparatoire zugelassen; diese ist jedoch eine weniger wissenschaftliche als didaktische Ausbildungsstätte, eingebettet in das Lyceum Louis-le-Grand. Galois bleibt an das ihn so bedrückende Umfeld gebunden. Als Student und Lehramtskandidat an der Ecole Préparatoire ist er immerhin automatisch Beamter des französischen Staates und verpflichtet sich mit seiner Einschreibung für sechs Jahre als Lehrer im Staatsdienst. 

Mathematisch entwickelt er sich in dieser Zeit dennoch weiter, wenn er auch mit erneuten Rückschlägen zurecht kommen muss: Cauchy verschlampt seine beiden Arbeiten, und eine weitere Arbeit, mit der Galois am \glqq Grand Prix de Mathématiques\grqq\ teilnehmen will und die er an Fourier, Sekretär an der Akademie und ebenfalls ein großer Mathematiker, übergibt, geht ebenfalls verloren: Fourier stirbt, und die Arbeit findet sich nicht mehr in seinem Nachlass. Hier häufen sich nun also die Rückschläge; und doch gibt es Hoffnung: So wird etwa eine Arbeit von Galois im renommierten \emph{Bulletin du Baron Férussac} veröffentlicht, neben Beiträgen von Cauchy, Poisson und Jacobi.  

Galois ist knappe 19 Jahre alt, als sich seine Lebens- mit der Zeitgeschichte zu vermengen beginnt. 1830, die Juli-Revolution: In Paris erheben sich Studenten und Arbeiter gegen Karl X., nachdem dieser sich anschickt, die Verfassung umzustoßen. Als Beamter darf Galois nicht an den Protesten gegen den König teilnehmen, eine unerträgliche Einschränkung für ihn, über die er sich vergeblich hinwegzusetzen versucht: Eine Flucht aus den verriegelten Toren der Schule gelingt nicht. Die Aufstände haben auch ohne Galois Erfolg: Im August löst Louis-Philippe von Orléans, der Bürgerkönig, Karl X. ab. 

Man darf aufgrund späterer Äußerungen vermuten, dass Galois dennoch enttäucht darüber ist, dass dies an ihm vorbeiläuft. Im Oktober 1830 jedenfalls tritt er dann der \glqq Societé des amis du peuple\grqq\ bei und gerät fortan immer häufiger in Konflikte mit dem Direktor Guigniault. Nach gewagten Äußerungen Galois' zum Thema der Uniformierung und Bewaffnung von Studenten -- Äußerungen, die auch in Form von Leserbriefen in verschiedenen Zeitungen erscheinen -- ist es genug: Guigniault verweist ihn im Dezember der Schule. Woraufhin sich Galois Hals über Kopf der 3. Batterie der Nationalgarde anschließt; seinem einengenden schulischen Umfeld endlich entflohen, das für die Weite seiner Ideen so wenig geeignet war, fließt nun ein enormer Teil an Energie in diese für ihn neue Ideenwelt. Auch in der Mathematik wird er zunehmend politischer; so verfasst er Leserbriefe, in denen er den gängigen mathematischen Lehrmethoden den Bankrott erklärt. Mathematik, so Galois, werde völlig falsch unterrichtet. Es könne nicht angehen, dass -- wie etwa bei den Sprachen -- Definitionen und Sätze nur heruntergebetet, auswendig gelernt und so nicht wirklich verstanden werden; es werde keinerlei Kreativität gefördert. Die Konsequenz dieser ketzerischen Ideen: endgültige Ausweisung aus dem Staatsdienst. Einer kurzen Episode als Privatlehrer ist kein Erfolg beschieden; zu unverständlich ist seine knappe Sprache noch immer den meisten anderen. 

Im Januar 1831 dann ein erneuter Versuch, seine Ideen zu publizieren; diesmal reicht Galois, der nach seinem Schulverweis Kurse an der Akademie besucht und dort durch zwar fundierte, aber sehr kritische Beiträge auffällt, die Arbeit bei Lacroix und Poisson ein. Er fügt ein kritisches Anschreiben hinzu, in dem er deutlich die Missstände bei der Akademie anprangert: sicher nicht die besten Voraussetzungen. Man wird Poisson zu gute halten müssen, dass er die Arbeit dennoch eingehend geprüft hat und die guten Ansätze erkannt haben muss. Dennoch erfolgt eine erneute Absage; jedoch mit der wohlwollenden Bemerkung, man werde, um sich ein echtes Bild machen zu können, die Veröffentlichung der gesamten Theorie Galois' abwarten müssen.  

Diese Absage seitens Poisson erhält Galois vermutlich, als er bereits im Gefängnis Sainte-Pélagie sitzt; dafür spricht das ebendort neu verfasste Vorwort zu der Arbeit, aus dem bitterste Ent\-täuschung spricht. Nach Sainte-Pélagie kommt Galois im Juli 1831, nachdem er zum zweiten Mal in diesem Jahr verhaftet worden ist. Die erste Verhaftung, bei einem durch die Societé finanzierten Bankett anlässlich des Freispruchs einiger Meuterer, erlebt Alexandre Dumas mit, der sich in seinen Memoiren daran erinnert, wie sich die Stimmung bei diesem Bankett immer mehr auflädt. Dumas spricht später von \glqq höchstgewagten Dialogen\grqq\ und \glqq außerordentlichen Toasts\grqq. Einen solchen Toast spricht auch Galois, unglücklicherweise mit einem Messer in der Hand und den Namen des Königs ausrufend: für die Gendarmerie eine Anstiftung zu einem Anschlag auf Louis-Philippe. Galois wird verhaftet, frei gesprochen, jedoch nach erneuter Verstrickung in politische Aktivitäten erneut verhaftet und schließlich nach Sainte-Pélagie gebracht. Dort erarbeitet er im Laufe eines halben Jahres weitere wesentliche Bestandteile seiner Theorie. Neben viel Zorn und Verbitterung findet sich in seinen Ausführungen auch ein äußerst kluges Wort: Zuzugeben, etwas nicht zu wissen, bedeute wirkliche Größe; niemals schade ein vermeintlich kluger Autor seinen Lesern so sehr, wie wenn er eine Schwierigkeit verbirgt. Ein Motto, das Galois beherzigt: Viele seiner Beiträge zeichnen sich gerade dadurch aus, dass er die \glqq richtigen Fragen\grqq\ stellt, auf die er zwar keine Antwort geben kann, die aber richtungsweisend sind. Gerade diese Art jedoch missfällt den Gutachtern oft, und man darf annehmen, dass die Ablehnung seiner Arbeiten zum großen Teil hierauf zurückzuführen ist. 

Galois' Leben, nun bereits auf der Zielgeraden, hält denn auch noch in letzter Minute eine amouröse Verstrickung bereit:  Stéphanie Poterin-Dumotel -- eine unglückliche Liebe, die außer zwei Brieffragmenten keinerlei Spuren hinterlässt. Es wird vermutet, dass er sie im Sanatorium Sieur Faultrier kennenlernt, wohin er aufgrund einer Choleraepidmie verlegt wird. Alles liegt hier völlig im Dunkeln, und sogar Stéphanies bloße Identität wurde erst in den 1960er Jahren durch eine eingehende Untersuchung der beiden Fragmente gelüftet. Entgegen vielen Gerüchten, in den ersten 150 Jahren nach Galois' Tod entstanden und wieder und wieder überliefert, geht man heute davon aus, dass es bei dem tragischen Duell nicht um Stéphanie geht. Ob der Hintergrund jedoch überhaupt privater oder vielleicht eher politischer Natur war, bleibt wohl für immer verborgen. 

Die Fakten: Am Morgen des 30. Mai 1832, einen Monat nach seiner Entlassung aus Sainte-Pélagie, kommt es zu einem Duell mit einem unbekannten Gegner; dabei wird Galois durch einen Schuss in den Unterleib so stark verletzt, dass er einen Tag später verstirbt. Den Beistand eines Priesters lehnt er noch ab. Ebenfalls eine Tatsache ist ein umfangreicher Brief an seinen Freund Auguste Chevalier, mit überwiegend mathematischem Inhalt, geschrieben fünf Tage vor dem Duell. Galois bezieht sich hier auf einige seiner bereits bekannten Ergebnisse, gibt Ergänzungen und weitere Beweisansätze. Im weiteren Sinne könnte man von einer Zusammenfassung seines Lebenswerkes sprechen. Indem er zeigen kann, dass die Fünfer-Permutationsguppe nicht auflösbar ist, erhält er einen eleganten Beweis von Abels Ergebnis und stellt es dadurch in einen viel größeren, erhellenderen Zusammenhang. Definitiv falsch ist die hartnäckige Legende, Galois habe all dies, seine komplette Theorie, in der Nacht vor dem Duell erstmals aufgeschrieben; das ist Romantik, vieles hatte er bereits gezeigt. Wahr ist allerdings, dass er in dem Brief seiner Hoffnung auf Nachhaltigkeit dieser Ergebnisse Ausdruck verleiht, ebenso wie er hofft, es werde sich jemand finden, der \glqq Erfolg beim Lesen all dieses Kuddelmuddels\grqq\ haben werde. 

Galois' Hoffnung wird erfüllt: Chevalier wird zum Nachlassverwalter und schickt Galois' Ergebnisse an viele Mathematiker. Nun endlich, posthum, wird das Genie Galois' erkannt; die in diesem Brief wie auch in seinen Arbeiten entwickelte Theorie trägt heute seinen Namen.  

\end{multicols}

\end{document}